\documentclass[a4paper,12pt,reqno]{amsart}
\usepackage{amsmath}
\usepackage{amssymb}
\usepackage{amsthm}

\usepackage{amscd}
\usepackage{amsfonts}
\usepackage{setspace}
\usepackage{version}

\title[Pickands' constant does not equal $1/\Gamma(1/\alpha)$, for small $\alpha$]{Pickands' constant $H_{\alpha}$ does not equal $1/\Gamma(1/\alpha)$, for small $\alpha$}
\author{Adam J Harper}
\address{Jesus College, Cambridge CB5 8BL, England}
\email{A.J.Harper@dpmms.cam.ac.uk}
\date{22nd April 2014}
\thanks{The author is supported by a research fellowship at Jesus College, Cambridge.}

\numberwithin{equation}{section}
\theoremstyle{plain}

\onehalfspace
\newcommand{\R}{\mathbb{R}}
\newcommand{\E}{\mathbb{E}}
\newcommand{\p}{\mathbb{P}}

\newtheorem{thmfromgp}{Theorem}
\newtheorem{piecelinthm}[thmfromgp]{Theorem}
\newtheorem{pickcorr}{Corollary}
\newtheorem{genthmfromgp}[thmfromgp]{Theorem}
\newtheorem{lem1}{Brownian Motion Lemma}
\newtheorem{lem2}[lem1]{Brownian Motion Lemma}
\newtheorem{lem3}[lem1]{Brownian Motion Lemma}
\newtheorem{prop1}{Proposition}
\newtheorem{prop2}[prop1]{Proposition}

\begin{document}

\maketitle

\begin{abstract}
Pickands' constants $H_{\alpha}$ appear in various classical limit results about tail probabilities of suprema of Gaussian processes. It is an often quoted conjecture that perhaps $H_{\alpha} = 1/\Gamma(1/\alpha)$ for all $0 < \alpha \leq 2$, but it is also frequently observed that this doesn't seem compatible with evidence coming from simulations.

We prove the conjecture is false for small $\alpha$, and in fact that $H_{\alpha} \geq (1.1527)^{1/\alpha}/\Gamma(1/\alpha)$ for all sufficiently small $\alpha$. The proof is a refinement of the ``conditioning and comparison'' approach to lower bounds for upper tail probabilities, developed in a previous paper of the author. Some calculations of hitting probabilities for Brownian motion are also involved.
\end{abstract}

% INTRODUCTION %%%%%%%%%%%%%%%%%%%%%%%%%%%%%
\section{Introduction}
In the author's paper~\cite{harpergp}, the following lower bound inequality was proved.
\begin{thmfromgp}[see Theorem 1 of Harper~\cite{harpergp}]
Let $n \geq 2$ and let $\{Z(t_{i})\}_{1 \leq i \leq n}$ be jointly multivariate normal random variables, each with mean zero and variance one. Suppose that the sequence is {\em stationary}, i.e. that $\E Z(t_{j})Z(t_{k}) = r(|j-k|)$ for some function $r$. Let $u \geq 1$, and suppose that:
\begin{itemize}
\item $r(m)$ is a decreasing non-negative function;

\item $r(1)(1+2u^{-2})$ is at most $1$.
\end{itemize}
Then $\p(\max_{1 \leq i \leq n} Z(t_{i}) > u)$ is
$$ \geq n \frac{e^{-u^{2}/2}}{40u} \min\left\{1,\sqrt{\frac{1-r(1)}{u^{2}r(1)}}\right\} \prod_{j=1}^{n-1} \Phi\left(u\sqrt{1-r(j)} \left(1+O\left(\frac{1}{u^{2}(1-r(j))}\right) \right) \right), $$
where $\Phi$ denotes the standard normal distribution function, and where the implicit constant in the ``big Oh'' notation is absolute (in particular, not depending on $\{Z(t_{i})\}_{1 \leq i \leq n}$), and could be found explicitly. 
\end{thmfromgp}
It turns out that Theorem 1 is almost sharp for some interesting collections of random variables $\{Z(t_{i})\}_{1 \leq i \leq n}$, for moderately sized $u$ (e.g. one can sometimes use Theorem 1 to identify $\E \max_{1 \leq i \leq n} Z(t_{i})$, up to second order terms). In the paper~\cite{harpergp}, Theorem 1 (or, more precisely, the ingredients of its proof) was used to obtain improved results in a probabilistic number theory problem. See the preprint~\cite{harpertypicalmax} for a (related) application to modelling the ``typical large values'' of the Riemann zeta function.

The proof of Theorem 1 breaks into two propositions. The first proposition was a {\em conditioning} step, in which $\p(\max_{1 \leq i \leq n} Z(t_{i}) > u)$ was lower bounded in terms of other probabilities involving conditioned versions of the $Z(t_{i})$. This was beneficial because, under the conditions on $r(m)$ imposed in Theorem 1, the correlation structure of the conditioned random variables could be lower bounded by a fairly nice correlation structure, corresponding to random variables constructed using random walks. The second proposition was a {\em comparison} step, in which Slepian's lemma was used to pass to random variables with the nicer lower bound correlation structure, and their behaviour was investigated using a simple result about the probability of Brownian motion remaining below a constant level.

In this paper we revisit the above argument, by requiring the Brownian motion in our comparison step to stay below a piecewise-linear function, rather than a constant. Most of this piecewise-linear function will be a negatively sloping line, which improves the bound by increasing the argument $u\sqrt{1-r(j)}$ in some of the product terms. Moreover, by choosing the height and slope of the function appropriately one can simultaneously improve the multiplier $\min\left\{1,\sqrt{\frac{1-r(1)}{u^{2}r(1)}}\right\}$. We will prove the following, slightly scary looking, result. In its statement, as well as in some of our later proofs, we use Vinogradov's notation $\gg$, meaning ``greater than, up to a multiplicative constant''. Thus a statement like $p(\alpha) \gg q(\alpha)$ means the same as $q(\alpha) = O(p(\alpha))$.
\begin{piecelinthm}
Let the situation be as in Theorem 1. In addition, let $C > 0$ and $K \geq 0$ and $1 \leq N \leq n-1$ be any parameters. Then $\p(\max_{1 \leq i \leq n} Z(t_{i}) > u)$ is
\begin{eqnarray}
& \gg & n \frac{e^{-u^{2}/2}}{u} \Phi\left(\frac{C/2 - Kr(N)/(1-r(N))}{\sqrt{r(N)/(1-r(N))}}\right) \min\left\{1,\frac{C(1-r(N))}{Kr(N)}\right\} \min\left\{1,\sqrt{\frac{C^{2}(1-r(1))}{r(1)}}\right\} \nonumber \\
&& \cdot \prod_{j=1}^{n-1} \Phi\left(u\sqrt{1-r(j)} \left(1 - \frac{C}{u} + \frac{K}{u} \min\{\frac{r(j)}{1-r(j)},\frac{r(N)}{1-r(N)}\} + O\left(\frac{1}{u^{2}(1-r(j))}\right) \right) \right) , \nonumber
\end{eqnarray}
where the implicit constants in the $\gg$ and ``big Oh'' notation are absolute (in particular, not depending on $\{Z(t_{i})\}_{1 \leq i \leq n}$), and could be found explicitly.
\end{piecelinthm}

Note that Theorem 1 follows from Theorem 2 by choosing $C=1/u$, $K = 0$, and $N=1$, say. As the reader will see later, the parameter $C$ in Theorem 2 may be thought of as a ``height'' parameter, the parameter $K$ may be thought of as a ``slope'' parameter, and the parameter $N$ may be thought of as a ``break'' parameter (where a boundary line of slope $-K$ changes into a horizontal line, of slope zero). Depending on the sizes of the $u\sqrt{1-r(j)}$ there may be other choices of the parameters that yield a stronger lower bound.

We shall prove Theorem 2 in $\S 2$ of this paper. For the benefit of a reader unfamiliar with the proof of Theorem 1, we will first prove a more general result (stated as Theorem 3, below) in which the conditioning and comparison steps are implemented, but the Brownian motion type term is left unanalysed. We then develop a few results about Brownian motion hitting probabilities for piecewise linear boundaries, and combine these with Theorem 3 to deduce Theorem 2. The proofs of fairly standard facts about Brownian motion are deferred to the appendix.

\vspace{12pt}
To illustrate the strength of Theorem 2, we turn to the Pickands constants application described in the title of this paper. Suppose that $\{Z(t)\}_{0 \leq t \leq h}$ is any mean zero, variance one, stationary Gaussian process indexed on the real line, whose covariance function $r(t) := \E Z(0)Z(t)$ satisfies
$$ r(t) = 1 - C|t|^{\alpha} + o(|t|^{\alpha}) \;\;\;\;\; \text{as} \; t \rightarrow 0 , $$
for some $C > 0$ and $0 < \alpha \leq 2$. An important theorem of Pickands~\cite{pickands1} asserts that, provided $\sup_{\epsilon \leq t \leq h} r(t) < 1$ for all $\epsilon > 0$, one has
$$ \lim_{u \rightarrow \infty} e^{u^{2}/2} u^{1-2/\alpha} \p\left(\sup_{0 \leq t \leq h} Z(t) > u \right) = \frac{hC^{1/\alpha} H_{\alpha}}{\sqrt{2\pi}} , $$
where $H_{\alpha}$ is the so-called {\em Pickands constant}.

It is a frequently quoted conjecture (see for example~\cite{burnmich}, noting that our $H_{\alpha}$ is written there as $\mathcal{H}_{\alpha/2}$) that perhaps $H_{\alpha} = 1/\Gamma(1/\alpha)$, and this is known to hold when $\alpha = 1, 2$, the only cases where the value of $H_{\alpha}$ is known exactly. But it is also frequently observed that, in general, this conjecture doesn't seem to match the behaviour predicted by simulations of random processes. In their preprint~\cite{diekyak}, Dieker and Yakir develop more practical Monte Carlo experiments for the investigation of $H_{\alpha}$, and state that ``...our simulation gives strong evidence that this conjecture is not correct... the confidence interval and [heuristic] error bounds are well above the curve [corresponding to $1/\Gamma(1/\alpha)$] for $\alpha$ in the range 1.6--1.8.'' D\c{e}bicki and Mandjes reproduce the conjecture in their open problems paper~\cite{debmand}, saying that it ``... lacks any firm heuristic support... [but] has not been falsified so far''.

Until recently, the best known lower bound for Pickands' constants (for small $\alpha$) was Michna's~\cite{michna} bound $H_{\alpha} \geq \frac{\alpha}{4\Gamma(1/\alpha)} (1/4)^{1/\alpha}$, which improved an earlier bound of D\c{e}bicki, Michna and Rolski~\cite{dmr} by a multiplicative factor of 2. In the author's paper~\cite{harpergp} this was improved using the techniques underlying Theorem 1, to show that for a small absolute constant $c > 0$ (which could be found explicitly) one has
$$ H_{\alpha} \geq \frac{c \alpha}{\Gamma(1/\alpha)} (1/2)^{1/\alpha} \;\;\;\;\; \forall 0 < \alpha \leq 2. $$
See the introduction to the author's paper~\cite{harpergp} for further references concerning bounds for Pickands' constants, and the introduction to Dieker and Yakir's preprint~\cite{diekyak} for further general references.

By applying Theorem 2 to suitable random variables $Z(t_{i})$, and making a good choice of the parameters $C,K,N$, we further improve the lower bound for $H_{\alpha}$ when $\alpha$ is small. This will be done in $\S 3$, below. In particular, we can show that the conjecture about $H_{\alpha}$ is false for small enough $\alpha$.
\begin{pickcorr}
There is an absolute constant $c > 0$, which could be found explicitly, such that
$$ H_{\alpha} \geq \frac{c \alpha^{5/2} (1.15279)^{1/\alpha}}{\Gamma(1/\alpha)} \;\;\;\;\; \forall 0 < \alpha \leq 2. $$
In particular, if $\alpha > 0$ is sufficiently small then $H_{\alpha} \geq (1.1527)^{1/\alpha}/\Gamma(1/\alpha)$.
\end{pickcorr}

\vspace{12pt}
The lower bound in Corollary 1 is essentially the best that seems to follow from Theorem 2, but is almost certainly not the best bound obtainable by our ``conditioning and comparison'' method. That is because Theorem 2 corresponds to the Brownian motion in our comparison step remaining below a negatively sloping line, and then a horizontal line, which is presumably not the best choice of boundary function for this application. When we apply Theorem 2 to prove Corollary 1, our choices of the parameters\footnote{More precisely, the choices of $K$ and $N$ are constrained by a few special sizes of $j$. It is easy to choose $C$ such that it isn't too small, but has a negligible effect in the product over $j$, which is more or less the best one can hope for since the product over $j$ is by far the hardest thing to control.} $C,K,N$ are dictated by the behaviour of $u\sqrt{1-r(j)}$ for a few special sizes of $j$, which suggests that if one considered a more complicated function one could work more carefully around those special ranges of $j$, and obtain a stronger bound. As a concrete (but probably quite fiddly) suggestion for further work, it would very likely lead to a stronger bound if one allowed a boundary function consisting of a negatively sloping line, and then another line with a different negative (rather than zero) slope, which would introduce an extra slope parameter into Theorem 2. In principle one can consider any boundary function, but estimating the relevant Brownian hitting probabilities may be impractical if the choice is too complicated.

The assumptions made on the correlation function $r(m)$ in Theorems 1 and 2 (and in the underlying ``conditioning and comparison'' arguments) are not very specialised, so the author believes there should also be several applications to other probabilistic problems. Some of these will be pursued in future work.

% SECTION 2 %%%%%%%%%%%%%%%%%%%%%%%%%%%%%
\section{Proof of Theorem 2}

\subsection{A more general result}
As mentioned in the introduction, to prove Theorem 2 we shall first state and prove a more general result, which encapsulates the conditioning and comparison arguments whilst leaving the Brownian motion (or, in fact, random walk) term for further analysis.
\begin{genthmfromgp}
Let the hypotheses be as in Theorem 1. Also let $(\delta(i))_{1 \leq i \leq n-1}$ be any real numbers. Then
\begin{eqnarray}
\p(\max_{1 \leq i \leq n} Z(t_{i}) > u) & \geq & n \frac{e^{-u^{2}/2}}{12u} \p\left(\sum_{j \leq i} \alpha_{j} Y_{j} \leq \delta(n-i)u \; \forall 1 \leq i \leq n-1 \right) \nonumber \\
&& \cdot \prod_{j=1}^{n-1} \Phi\left(u\sqrt{1-r(j)} \left(1 - \delta(j) + O\left(\frac{1}{u^{2}(1-r(j))}\right) \right) \right) , \nonumber
\end{eqnarray}
where the $Y_{j}$ are independent standard normal random variables, and the $\alpha_{j}$ defined by
$$ \sum_{j \leq i} \alpha_{j}^{2} := \frac{r(n-i)}{1-r(n-i)} . $$
The implicit constant  in the ``big Oh'' notation is absolute (in particular, not depending on $\{Z(t_{i})\}_{1 \leq i \leq n}$), and could be found explicitly.
\end{genthmfromgp}

Notice that $r(n-i)/(1-r(n-i))$ is an increasing function of $1 \leq i \leq n-1$, since $r(m)$ is assumed to be a decreasing function.

\vspace{12pt}
Theorem 3 may be extracted from the proofs of Propositions 1 and 2 in the author's paper~\cite{harpergp}, but we shall recap the main details.

Since we assume the $\{Z(t_{i})\}_{1 \leq i \leq n}$ are stationary, we see $\p(\max_{1 \leq i \leq n} Z(t_{i}) > u)$ is
\begin{eqnarray}
& = & \sum_{m=1}^{n} \p(Z(t_{m}) > u, Z(t_{j}) \leq u \; \forall 1 \leq j \leq m-1) \nonumber \\
& \geq & n \p(Z(t_{n}) > u, Z(t_{j}) \leq u \; \forall 1 \leq j \leq n-1) \nonumber \\
& = & n \p\left(Z(t_{n}) > u, \frac{Z(t_{j}) - r(n-j)Z(t_{n})}{\sqrt{1-r(n-j)^{2}}} \leq \frac{u - r(n-j)Z(t_{n})}{\sqrt{1-r(n-j)^{2}}} \; \forall 1 \leq j \leq n-1 \right) , \nonumber
\end{eqnarray}
and it is easy to check that the random variables $V_{j} := \frac{Z(t_{j}) - r(n-j)Z(t_{n})}{\sqrt{1-r(n-j)^{2}}}$ satisfy
$$ \E V_{j} = 0, \;\;\;\;\; \E V_{j}^{2}=1, \;\;\;\;\; \E V_{j}V_{k} = \frac{r(|j-k|) - r(n-j)r(n-k)}{\sqrt{1-r(n-j)^{2}} \sqrt{1-r(n-k)^{2}}}, \;\;\;\;\; \E V_{j} Z(t_{n}) = 0 . $$
In particular, since the $\{Z(t_{i})\}_{1 \leq i \leq n}$ were assumed to be jointly normal and since $\E V_{j} Z(t_{n}) = 0$ we know the $V_{j}$ are all independent of $Z(t_{n})$, so {\em conditioning} shows
\begin{eqnarray}
\p(\max_{1 \leq i \leq n} Z(t_{i}) > u) & \geq & n \int_{u}^{u+1/u} \p(V_{j} \leq \frac{u - r(n-j)x}{\sqrt{1-r(n-j)^{2}}} \; \forall 1 \leq j \leq n-1) \frac{e^{-x^{2}/2}}{\sqrt{2\pi}} dx \nonumber \\
& \geq & n \frac{e^{-(u+1/u)^{2}/2}}{u\sqrt{2\pi}} \inf_{u \leq x \leq u+1/u} \p(V_{j} \leq \frac{u - r(n-j)x}{\sqrt{1-r(n-j)^{2}}} \; \forall 1 \leq j \leq n-1) . \nonumber
\end{eqnarray}
Since we assume that $r(m)$ is a non-negative function, the infimum is attained when $x=u+1/u$, so a quick calculation (using our assumption that $u \geq 1$) yields the simplified lower bound
$$ \p(\max_{1 \leq i \leq n} Z(t_{i}) > u) \geq n \frac{e^{-u^{2}/2}}{12u} \p\left(V_{j} \leq \frac{u(1 - r(n-j)(1+u^{-2}))}{\sqrt{1-r(n-j)^{2}}} \; \forall 1 \leq j \leq n-1 \right) . $$

(The above corresponds to choosing $H=1/u$ in Proposition 1 in the author's paper~\cite{harpergp}, and using stationarity to slightly simplify the form of the bound.)

Next, observe that
$$ \E V_{j} V_{k} \geq \frac{r(n-\min\{j,k\}) (1- r(n-\max\{j,k\}))}{\sqrt{1-r(n-j)^{2}} \sqrt{1-r(n-k)^{2}}} \;\;\;\;\; \forall 1 \leq j,k \leq n-1, $$
since $r(|j-k|) \geq r(n-\min\{j,k\})$ (as $r(m)$ is assumed to be a decreasing function). Therefore by Slepian's Lemma (see e.g. Comparison Inequality 2 in the author's paper~\cite{harpergp}), if $\{X_{j}\}_{1 \leq j \leq n-1}$ are mean zero, variance one, jointly normal random variables such that $\E X_{j}X_{k} = \frac{r(n-\min\{j,k\}) (1- r(n-\max\{j,k\}))}{\sqrt{1-r(n-j)^{2}} \sqrt{1-r(n-k)^{2}}}$ for all $j \neq k$, then we have the {\em comparison} lower bound
$$ \p\left(V_{j} \leq \frac{u(1 - r(n-j)(1+u^{-2}))}{\sqrt{1-r(n-j)^{2}}} \forall j \leq n-1 \right) \geq \p\left(X_{j} \leq \frac{u(1 - r(n-j)(1+u^{-2}))}{\sqrt{1-r(n-j)^{2}}} \forall j \leq n-1 \right) . $$
In the proof of Proposition 2 in the author's paper~\cite{harpergp} (with the choices $c_{j} = r(n-j)$ and $d_{j} = 1 - r(n-j)$), it is shown by construction that such random variables $X_{j}$ always exist, and that $\p\left(X_{j} \leq \frac{u(1 - r(n-j)(1+u^{-2}))}{\sqrt{1-r(n-j)^{2}}} \forall j \leq n-1 \right)$ is
$$ = \p\left(Z_{i} \leq \frac{u(1 - r(n-i)(1+u^{-2})) - (1-r(n-i))\sum_{j \leq i} \alpha_{j}Y_{j}}{\sqrt{1-r(n-i)}} \; \forall i \leq n-1 \right) , $$
where the $Z_{i}$ and the $Y_{j}$ are all {\em independent} standard normal random variables, and the real numbers $\alpha_{j}$ are defined as in Theorem 3.

Finally, using independence, for any real $(\delta(i))_{1 \leq i \leq n-1}$ the above probability is
\begin{eqnarray}
& \geq & \p\left(\sum_{j \leq i} \alpha_{j} Y_{j} \leq \delta(n-i)u \; \forall 1 \leq i \leq n-1 \right) \nonumber \\
&& \cdot \p\left(Z_{i} \leq \frac{u((1 - r(n-i))(1-\delta(n-i)) - r(n-i)u^{-2})}{\sqrt{1-r(n-i)}} \; \forall i \leq n-1 \right) \nonumber \\
& = & \p\left(\sum_{j \leq i} \alpha_{j} Y_{j} \leq \delta(n-i)u \; \forall 1 \leq i \leq n-1 \right) \nonumber \\
&& \cdot \prod_{i=1}^{n-1} \Phi\left( u\sqrt{1 - r(n-i)}\left(1 - \delta(n-i) - \frac{r(n-i)}{u^{2}(1-r(n-i))}\right) \right) , \nonumber
\end{eqnarray}
from which Theorem 3 follows.
\begin{flushright}
Q.E.D.
\end{flushright}

\subsection{Some calculations with Brownian motion}
In this subsection we shall perform a few calculations involving Brownian motion, which will ultimately supply a lower bound for the term $\p\left(\sum_{j \leq i} \alpha_{j} Y_{j} \leq \delta(n-i)u \; \forall 1 \leq i \leq n-1 \right)$ in Theorem 3 (for a special choice of the numbers $(\delta(i))_{1 \leq i \leq n-1}$).

We begin by stating two lower bounds for the probability of Brownian motion remaining below a negatively sloping line segment.
\begin{lem1}
Let $a > 0$, $b < 0$, and $t > 0$. Suppose that $|b\sqrt{t}|$ is sufficiently large. Then if $\{W_{s}\}_{s \geq 0}$ denotes a standard Brownian motion (started from zero), we have
$$ \p(W_{s} \leq a+bs \; \forall 0 \leq s \leq t) \gg \min\left\{1,\frac{a}{|bt|}\right\} \Phi(\frac{a+bt}{\sqrt{t}}) , $$
where the constant implicit in the $\gg$ notation is absolute.
\end{lem1}

\begin{lem2}
Let $H > 0$ be any fixed constant. Let $a > 0$, $b \leq 0$, and $t > 0$, and suppose that $|b\sqrt{t}| \leq H$. Then if $\{W_{s}\}_{s \geq 0}$ denotes a standard Brownian motion (started from zero), we have
$$ \p(W_{s} \leq a+bs \; \forall 0 \leq s \leq t) \gg_{H} \min\left\{1,\frac{a}{\sqrt{t}}\right\} , $$
where the constant implicit in the $\gg_{H}$ notation depends on $H$ only.
\end{lem2}

Brownian Motion Lemmas 1 and 2 are consequences of the well known explicit formula for hitting probabilities of a sloping line by Brownian motion, together with a little analysis to simplify the resulting expressions. For the sake of completeness, proofs of these lemmas are included in the appendix.

Combining Brownian Motion Lemmas 1 and 2, we can deduce the following result.
\begin{lem3}
Let $a > 0$, $b \leq 0$, and $0 < t_{0} < t$. Then if $\{W_{s}\}_{s \geq 0}$ denotes a standard Brownian motion (started from zero), we have
$$ \p(W_{s} \leq a+b\min\{s,t_{0}\} \; \forall 0 \leq s \leq t) \gg \min\left\{1,\frac{a}{|bt_{0}|}\right\} \Phi(\frac{a/2+bt_{0}}{\sqrt{t_{0}}}) \min\left\{1,\frac{a}{\sqrt{t}}\right\} , $$
where the constant implicit in the $\gg$ notation is absolute.
\end{lem3}
To prove Brownian Motion Lemma 3 we distinguish two cases. Let $H > 1$ be a sufficiently large constant that Brownian Motion Lemma 1 is applicable when $|b\sqrt{t_{0}}| \geq H$. Then if $|b\sqrt{t_{0}}| \geq H$, we observe that
\begin{eqnarray}
&& \p(W_{s} \leq a+b\min\{s,t_{0}\} \; \forall 0 \leq s \leq t) \nonumber \\
& \geq & \p(W_{s} \leq a/2+bs \; \forall 0 \leq s \leq t_{0}, \;\;\; \text{and} \;\;\; (W_{s}-W_{t_{0}}) \leq a/2 \; \forall t_{0} < s \leq t) \nonumber \\
&  = & \p(W_{s} \leq a/2+bs \; \forall 0 \leq s \leq t_{0}) \cdot \p(B_{s} \leq a/2 \; \forall 0 \leq s \leq t-t_{0}) \nonumber \\
& \gg & \min\left\{1,\frac{a}{|bt_{0}|}\right\} \Phi(\frac{a/2+bt_{0}}{\sqrt{t_{0}}}) \min\left\{1,\frac{a}{\sqrt{t-t_{0}}}\right\} , \nonumber
\end{eqnarray}
where $B_{s}$ denotes another standard Brownian motion. Here the final inequality uses Brownian Motion Lemma 1, and then the well known fact that $\max_{0 \leq s \leq t-t_{0}} B_{s} \sim |N(0,t-t_{0})|$ (or, alternatively, Brownian Motion Lemma 2 with $b=0$).

The other case is where $|b\sqrt{t_{0}}| < H$. Let $\tilde{b} := \min\{b,-1/\sqrt{t_{0}}\}$. Then using Brownian Motion Lemma 2, we have (remembering that $b, \tilde{b}$ are non-positive)
\begin{eqnarray}
&& \p(W_{s} \leq a+b\min\{s,t_{0}\} \; \forall 0 \leq s \leq t) \nonumber \\
& \geq & \p(W_{s} \leq a/2+2\tilde{b}s \; \forall 0 \leq s \leq t_{0}, \;\;\; \text{and} \;\;\; (W_{s}-W_{t_{0}}) \leq a/2 + |\tilde{b}|t_{0} \; \forall t_{0} < s \leq t) \nonumber \\
&  = & \p(W_{s} \leq a/2+2\tilde{b}s \; \forall 0 \leq s \leq t_{0}) \cdot \p(B_{s} \leq a/2 + |\tilde{b}|t_{0} \; \forall 0 \leq s \leq t-t_{0}) \nonumber \\
& \gg_{H} & \min\left\{1,\frac{a}{\sqrt{t_{0}}}\right\} \min\left\{1,\frac{a/2 + |\tilde{b}|t_{0}}{\sqrt{t-t_{0}}}\right\} \nonumber \\
& \gg & \min\left\{1,\frac{a}{\sqrt{t_{0}}}\right\} \min\left\{1,\frac{a + \sqrt{t_{0}}}{\sqrt{t-t_{0}}}\right\} \gg \min\left\{1,\frac{a}{\sqrt{t}}\right\} . \nonumber
\end{eqnarray}
Here the final inequality follows by considering whether $a \geq \sqrt{t_{0}}$ or not. We also observe that, since $H$ is now an absolute constant (determined only by the meaning of ``sufficiently large'' in the statement of Brownian Motion Lemma 1), we can drop the subscript on the $\gg_{H}$ notation that denotes dependence on $H$. 

To summarise, in both cases we have shown, as claimed, that
$$ \p(W_{s} \leq a+b\min\{s,t_{0}\} \; \forall 0 \leq s \leq t) \gg \min\left\{1,\frac{a}{|bt_{0}|}\right\} \Phi(\frac{a/2+bt_{0}}{\sqrt{t_{0}}}) \min\left\{1,\frac{a}{\sqrt{t}}\right\} . $$
\begin{flushright}
Q.E.D.
\end{flushright}

We conclude this subsection with two (fairly obvious) remarks.

Firstly, the proof of Brownian Motion Lemma 3 is a bit wasteful on certain ranges of the parameters $a,b,t_{0},t$. However, a sharper bound would be more complicated to state, and (it seems) of little additional use for the ultimate Pickands constants application.

Secondly, the main reason for examining linear and piecewise linear boundaries here (which translates into Theorem 2, as will be seen in the next subsection) is simply that it is fairly easy to work with them, because of the corresponding explicit formula for Brownian motion hitting probabilities. It is quite reasonable to think that, in any given application, the best choice of the numbers $(\delta(i))_{1 \leq i \leq n-1}$ in Theorem 3 will not correspond to a piecewise linear boundary, although Theorem 2 (and even Theorem 1) do seem to perform quite well in various applications.

\subsection{Putting everything together}
The proof of Theorem 2 is completed by combining Theorem 3 with Brownian Motion Lemma 3. Indeed, if we simply choose
$$ \delta(i) := \frac{C}{u} - \frac{K}{u} \min\left\{\frac{r(i)}{1-r(i)}, \frac{r(N)}{1-r(N)}\right\}  $$
in Theorem 3 then the product over $j$ there is as required for Theorem 2, whilst
\begin{eqnarray}
&& \p\left(\sum_{j \leq i} \alpha_{j} Y_{j} \leq \delta(n-i)u \; \forall 1 \leq i \leq n-1 \right) \nonumber \\
& = & \p\left(\sum_{j \leq i} \alpha_{j} Y_{j} \leq C - K\min\left\{\sum_{j \leq i} \alpha_{j}^{2}, \frac{r(N)}{1-r(N)} \right\}  \; \forall 1 \leq i \leq n-1 \right) \nonumber \\
& \geq & \p\left(W_{s} \leq C - K\min\left\{s,\frac{r(N)}{1-r(N)}\right\} \; \forall 0 \leq s \leq \frac{r(1)}{1-r(1)} \right), \nonumber
\end{eqnarray}
since $\sum_{j \leq i} \alpha_{j}^{2} = \frac{r(n-i)}{1-r(n-i)}$ and since we always have $(\sum_{j \leq i} \alpha_{j} Y_{j})_{1 \leq i \leq n-1} \stackrel{d}{=} (W_{\sum_{j \leq i} \alpha_{j}^{2}})_{1 \leq i \leq n-1}$ (where $\stackrel{d}{=}$ denotes equality in distribution). Brownian Motion Lemma 3 then shows this probability is
$$ \gg \min\left\{1,\frac{C (1-r(N))}{Kr(N)}\right\} \Phi(\frac{C/2 - Kr(N)/(1-r(N))}{\sqrt{r(N)/(1-r(N))}}) \min\left\{1, \sqrt{\frac{C^{2} (1-r(1))}{r(1)}}\right\} , $$
as required for Theorem 2.
\begin{flushright}
Q.E.D.
\end{flushright}

% SECTION 3 %%%%%%%%%%%%%%%%%%%%%%%%%%%%%
\section{Proof of Corollary 1}

\subsection{Overview of the argument}
It doesn't really require any further ideas to deduce Corollary 1 from Theorem 2, but the details of the calculation are quite involved. To try to clarify things, in this subsection we describe the collections of random variables to which Theorem 2 will be applied, and divide the task of deducing Corollary 1 into two further propositions. Those propositions will be proved in the following subsections.

Let $0 < \alpha < 2$, and let $\{Z(t)\}_{t \geq 0}$ be a mean zero, variance one, stationary Gaussian process with covariance function
$$ r(t) := \E Z(0)Z(t) = \frac{1}{2}\left(e^{\alpha t/2} + e^{-\alpha t/2} - (e^{t/2}-e^{-t/2})^{\alpha}\right), \;\;\; t \geq 0. $$
Such Gaussian processes were constructed by Shao~\cite{shao} in his work on Pickands' constants, by suitably reparametrising fractional Brownian motion. It is easy to check that
$$ r(t) = 1 - t^{\alpha}/2 + O(t^{2}) \;\;\;\;\; \text{as} \; t \rightarrow 0 , $$
and therefore by Pickands' theorem (as stated in the introduction) we have
$$ H_{\alpha} = 2^{1/\alpha} \sqrt{2\pi} \lim_{u \rightarrow \infty} e^{u^{2}/2} u^{1-2/\alpha} \p\left(\sup_{0 \leq t \leq 1} Z(t) > u \right) . $$

Let us make two further remarks. In our proofs it will be convenient to assume that $\alpha < \alpha_{0}$, for a certain small number $\alpha_{0} > 0$. For any {\em fixed} value $\alpha_{0} > 0$, Corollary 1 holds for all $\alpha_{0} \leq \alpha \leq 2$ as a consequence of the existing lower bounds for $H_{\alpha}$ (e.g. the bound due to Michna~\cite{michna}), provided the constant $c > 0$ in Corollary 1 is small enough. Thus we can indeed restrict our arguments to the case $\alpha < \alpha_{0}$, where $\alpha_{0}$ is a small fixed constant. (An explicit permissible choice of $\alpha_{0}$ could be found by working very carefully through all our proofs. The author believes that setting $\alpha_{0}=1/400$ is more than sufficient, but hasn't checked this fully since it doesn't affect the overall shape of our bounds.) Let us also recall that, by Stirling's formula, $\Gamma(1/\alpha) \sim \sqrt{2\pi} (1/\alpha)^{1/\alpha - 1/2} e^{-1/\alpha}$ as $\alpha \rightarrow 0$. So in order to prove Corollary 1, it will suffice to prove the following result.
\begin{prop1}
There exists a small constant $\alpha_{0} > 0$, which could be found explicitly, such that the following is true.

Let $0 < \alpha \leq \alpha_{0}$, and let $\{Z(t)\}_{t \geq 0}$ be the Gaussian process described above. Then provided $u$ is sufficiently large in terms of $\alpha$,
$$ \p\left(\sup_{0 \leq t \leq 1} Z(t) > u \right) \gg \frac{e^{-u^{2}/2}}{u} u^{2/\alpha} \frac{1}{2^{1/\alpha}} \alpha^{2} (1.15279 e \alpha)^{1/\alpha} , $$
where the constant implicit in the $\gg$ notation is absolute.
\end{prop1}

We shall ultimately use Theorem 2 to prove Proposition 1. To do this, note first that $r(t)$ is a decreasing non-negative function of $t \geq 0$, which is easily checked by calculating $r'(t)$ (as was done at the beginning of $\S 5$ in the author's paper~\cite{harpergp}.) Next, for any integer $M=M(u,\alpha) \geq 1$ we obviously have
$$ \p\left(\sup_{0 \leq t \leq 1} Z(t) > u \right) \geq \p\left(\max_{1 \leq i \leq M} Z(i/M) > u \right) , $$
and the random variables $\{Z(i/M)\}_{1 \leq i \leq M}$ will satisfy all the conditions of Theorem 2 provided that
$$ r(1/M)(1+2u^{-2}) \leq 1. $$
Let us choose $M = \lfloor (bu^{2}\alpha/2)^{1/\alpha} \rfloor$, where $\lfloor \cdot \rfloor$ denotes integer part and where $1 \leq b \leq 100$ is a constant whose optimal value (from the point of view of proving Proposition 1) will be determined later. If $u$ is sufficiently large in terms of $\alpha$ we can make $M$ arbitrarily large, and so (since $\alpha \leq \alpha_{0} \leq 1/400 \leq 1/4b$) we see $r(1/M)(1+2u^{-2})$ is
\begin{eqnarray}
= \left(1-\frac{1}{2M^{\alpha}} + O\left(\frac{1}{M^{2}}\right)\right) (1+2u^{-2}) \leq (1-\frac{1}{4M^{\alpha}} ) (1+2u^{-2}) & \leq & (1-\frac{1}{2bu^{2}\alpha} ) (1+2u^{-2}) \nonumber \\
& < & 1 . \nonumber
\end{eqnarray}
Thus the random variables $\{Z(i/M)\}_{1 \leq i \leq M}$ do satisfy all the conditions of Theorem 2, provided $u$ is sufficiently large in terms of $\alpha$.

We must still decide how to choose the ``height'', ``slope'' and ``break'' parameters $C > 0$, $K \geq 0$ and $1 \leq N \leq M-1$ in Theorem 2, in order to obtain the best lower bound we can. We will divide the task of doing this into two parts: in the next subsection we shall prove the following proposition, in the course of which we will choose $C$ and also the rough forms of $K$ and $N$ (in terms of two further parameters $\kappa, Y$); afterwards we will finetune the choices of $\kappa$ and $Y$, and also make the best choice of $b$ that we can, to finally deduce Proposition 1.
\begin{prop2}
There exists a small constant $\alpha_{0} > 0$, which could be found explicitly, such that the following is true.

Let $0 < \alpha \leq \alpha_{0}$, let $\{Z(t)\}_{t \geq 0}$ be the Gaussian process described above, let $1 \leq b \leq 100$, and set $M = \lfloor (bu^{2}\alpha/2)^{1/\alpha} \rfloor$. Finally, let $1/1000 \leq \kappa \leq 1000$ and $1 \leq Y \leq 1000$ be any parameters. Then provided $u$ is sufficiently large in terms of $\alpha$,
\begin{eqnarray}
\p\left(\max_{1 \leq i \leq M} Z(i/M) > u \right) & \gg & M \frac{e^{-u^{2}/2}}{u} \alpha^{3/2} \Phi\left(-(1+O(\alpha))\kappa \sqrt{\frac{b}{\alpha Y}} \right) \nonumber \\
&& \cdot \prod_{j \leq M^{1/4}} \Phi\left((1+O(\alpha))\sqrt{\frac{j^{\alpha}}{b\alpha}} \left(1 + \kappa b  \min\left\{\frac{1}{j^{\alpha}},\frac{1}{Y}\right\} \right) \right) , \nonumber
\end{eqnarray}
where the constants implicit in the $\gg$ and ``big Oh'' notations are absolute.
\end{prop2}

\subsection{Proof of Proposition 2}
Proposition 2 is a messy but straightforward deduction from Theorem 2, repeatedly using the fact that our underlying covariance function satisfies
$$ r(t) = 1 - t^{\alpha}/2 + O(t^{2}) \;\;\;\;\; \text{as} \; t \rightarrow 0 . $$
Many details of the deduction are the same as in $\S 5$ of the author's paper~\cite{harpergp}.

It is helpful first to observe that if $j > M^{1/4}$ we have
$$ u\sqrt{1-r(j/M)} \geq u\sqrt{1-r(M^{-3/4})} = u\sqrt{M^{-3\alpha/4}/2 + O(M^{-3/2})} \geq \frac{u}{2M^{3\alpha/8}} \geq u^{1/4}, $$
say, since $r(t)$ is decreasing and $M = \lfloor (bu^{2}\alpha/2)^{1/\alpha} \rfloor$ is large (and $\alpha \leq \alpha_{0} \leq 1/4b$ is small). Therefore for {\em any} choice of $0 < C \leq u/2$ (note the very weak upper bound restriction on $C$) and $K \geq 0$ and $1 \leq N \leq M-1$, we will have
\begin{eqnarray}
&& \prod_{M^{1/4} < j \leq M} \Phi\left(u\sqrt{1-r(\frac{j}{M})} \left(1 - \frac{C}{u} + \frac{K}{u} \min\{\frac{r(\frac{j}{M})}{1-r(\frac{j}{M})},\frac{r(\frac{N}{M})}{1-r(\frac{N}{M})}\} + O\left(\frac{1}{u^{2}(1-r(\frac{j}{M}))}\right) \right) \right) \nonumber \\
& \geq & \prod_{M^{1/4} < j \leq M} \Phi\left(u\sqrt{1-r(\frac{j}{M})} \left(1/2 + O\left(\frac{1}{u^{1/2}}\right) \right) \right) \geq \left(\Phi(u^{1/4}/4) \right)^{M} , \nonumber
\end{eqnarray}
since $u$ is large. Since we have $M \leq (50u^{2}\alpha)^{1/\alpha}$, (a power of $u$, for any fixed $\alpha$), we will have $\left(\Phi(u^{1/4}/4) \right)^{M} \geq \left(1-e^{-u^{1/2}/32}\right)^{M} \geq 1/2$, say, provided $u$ is large enough in terms of $\alpha$. So in Theorem 2 we only need to deal with the part of the product where $j \leq M^{1/4}$.

Next, when $j \leq M^{1/4}$ we have
\begin{eqnarray}
u\sqrt{1-r(j/M)} = u\sqrt{\frac{j^{\alpha}}{2M^{\alpha}} + O(\frac{j^{2}}{M^{2}})} & = & \left(1+O\left(\frac{j^{2-\alpha}}{M^{2-\alpha}}\right) \right) u\sqrt{\frac{j^{\alpha}}{2M^{\alpha}}} \nonumber \\
& = & (1+O(M^{-3/4})) u\sqrt{\frac{j^{\alpha}}{2M^{\alpha}}}, \nonumber
\end{eqnarray}
say, since $\alpha$ is small. The ``big Oh'' term here is much better than we need, and for convenience of writing later we take a very crude approach and note it is certainly $O(\alpha/j^{\alpha})$, provided $u$ is large enough in terms of $\alpha$. So we have
$$ u\sqrt{1-r(j/M)} = (1+O(\alpha/j^{\alpha})) u\sqrt{\frac{j^{\alpha}}{2M^{\alpha}}} \geq (1+O(\alpha/j^{\alpha})) u\sqrt{\frac{j^{\alpha}}{bu^{2}\alpha}} = (1+O(\alpha/j^{\alpha})) \sqrt{\frac{j^{\alpha}}{b\alpha}}  , $$ 
by definition of $M$. Inserting all this in Theorem 2, we see $\p\left(\max_{1 \leq i \leq M} Z(i/M) > u \right)$ is
\begin{eqnarray}
& \gg & M \frac{e^{-u^{2}/2}}{u} \Phi\left(\frac{C/2 - Kr(\frac{N}{M})/(1-r(\frac{N}{M}))}{\sqrt{r(\frac{N}{M})/(1-r(\frac{N}{M}))}}\right) \min\left\{1,\frac{C(1-r(\frac{N}{M}))}{Kr(\frac{N}{M})}\right\} \min\left\{1,\sqrt{\frac{C^{2}(1-r(\frac{1}{M}))}{r(\frac{1}{M})}}\right\} \nonumber \\
&& \cdot \prod_{j \leq M^{1/4}} \Phi\left(\left(1+O\left(\frac{\alpha}{j^{\alpha}}\right) \right)\sqrt{\frac{j^{\alpha}}{b\alpha}} \left(1 - \frac{C}{u} + \frac{K}{u} \min\{\frac{r(\frac{j}{M})}{1-r(\frac{j}{M})},\frac{r(\frac{N}{M})}{1-r(\frac{N}{M})}\} + O\left(\frac{\alpha}{j^{\alpha}}\right) \right) \right) \nonumber \\
& \gg & M \frac{e^{-u^{2}/2}}{u} \Phi\left(\frac{C/2 - Kr(\frac{N}{M})/(1-r(\frac{N}{M}))}{\sqrt{r(\frac{N}{M})/(1-r(\frac{N}{M}))}}\right) \min\left\{1,\frac{C(1-r(\frac{N}{M}))}{Kr(\frac{N}{M})}\right\} \min\left\{1,\frac{C}{u} \sqrt{\frac{1}{b\alpha}}\right\} \nonumber \\
&& \cdot \prod_{j \leq M^{1/4}} \Phi\left(\left(1+O\left(\frac{\alpha}{j^{\alpha}}\right) \right)\sqrt{\frac{j^{\alpha}}{b\alpha}} \left(1 - \frac{C}{u} + \frac{K}{u} \min\{\frac{r(\frac{j}{M})}{1-r(\frac{j}{M})},\frac{r(\frac{N}{M})}{1-r(\frac{N}{M})}\} + O\left(\frac{\alpha}{j^{\alpha}}\right) \right) \right) , \nonumber
\end{eqnarray}
where $0 < C \leq u/2$, $K \geq 0$ and $1 \leq N \leq M-1$ are still to be chosen.

To get a rough idea of how we should select our parameters, note that if $C,K \approx 0$ then the product over $j$ looks roughly like $\prod_{j \leq M^{1/4}}\Phi(\sqrt{j^{\alpha}/b\alpha})$, which is $\approx 1$ provided $b \leq e/2$, is very small if $b > e/2$, and moreover is dominated by those terms $j \leq (1000b)^{1/\alpha}$, say. (See $\S 5$ of the author's paper~\cite{harpergp} for an analysis of the behaviour of the product. We will also analyse it extensively in the next subsection.) So if we want to choose $b$ larger, as we do to prove Proposition 1, we need to choose $K$ such that $(K/u)r(N/M)/(1-r(N/M))$ is at least a large constant. Assuming that we shall choose $N \leq M^{1/4}$, (which seems sensible both to increase the size of $r(N/M)/(1-r(N/M))$, and because the size of the product is mostly determined by small $j$), our previous calculations show that
\begin{eqnarray}
\frac{K}{u} \min\{\frac{r(\frac{j}{M})}{1-r(\frac{j}{M})},\frac{r(\frac{N}{M})}{1-r(\frac{N}{M})}\} & = & Ku\min\{\frac{r(\frac{j}{M})}{u^{2}(1-r(\frac{j}{M}))},\frac{r(\frac{N}{M})}{u^{2}(1-r(\frac{N}{M}))}\} \nonumber \\
& = & Ku\min\{\frac{r(\frac{j}{M})}{(1+O(\alpha/j^{\alpha}))j^{\alpha}/(b\alpha)},\frac{r(\frac{N}{M})}{(1+O(\alpha/N^{\alpha}))N^{\alpha}/(b\alpha)}\} \nonumber \\
& = & Ku (b\alpha) \min\{\frac{1+O(\alpha/j^{\alpha})}{j^{\alpha}},\frac{1+O(\alpha/N^{\alpha})}{N^{\alpha}}\} . \nonumber
\end{eqnarray}
Here the final equality uses the fact that $r(j/M) = 1 + O((j/M)^{\alpha}) = 1 + O(\alpha/j^{\alpha})$ when $j \leq M^{1/4}$ and $u$ is large.

Motivated by all of this, let us take $K = \kappa/(u\alpha)$ and $N = Y^{1/\alpha}$, where $1/1000 \leq \kappa \leq 1000$ and $1 \leq Y \leq 1000$, say. Let us also set $C=u\alpha$, which certainly satisfies our earlier restriction that $0 < C \leq u/2$. With these choices we see
$$ \frac{K}{u} \min\{\frac{r(\frac{j}{M})}{1-r(\frac{j}{M})},\frac{r(\frac{N}{M})}{1-r(\frac{N}{M})}\} = (1+O(\alpha)) \kappa b \min\left\{\frac{1}{j^{\alpha}},\frac{1}{N^{\alpha}}\right\} = (1+O(\alpha)) \kappa b \min\left\{\frac{1}{j^{\alpha}},\frac{1}{Y}\right\}  , $$
and so Theorem 2 implies, as above, that $\p\left(\max_{1 \leq i \leq M} Z(i/M) > u \right)$ is
\begin{eqnarray}
& \gg & M \frac{e^{-u^{2}/2}}{u} \Phi\left(\frac{\alpha/2 - (K/u)r(\frac{N}{M})/(1-r(\frac{N}{M}))}{\sqrt{r(\frac{N}{M})/u^{2}(1-r(\frac{N}{M}))}}\right) \min\left\{1,\frac{u\alpha(1-r(\frac{N}{M}))}{Kr(\frac{N}{M})}\right\} \sqrt{\frac{\alpha}{b}} \nonumber \\
&& \cdot \prod_{j \leq M^{1/4}} \Phi\left(\left(1+O\left(\frac{\alpha}{j^{\alpha}}\right) \right)\sqrt{\frac{j^{\alpha}}{b\alpha}} \left(1 + \frac{K}{u} \min\{\frac{r(\frac{j}{M})}{1-r(\frac{j}{M})},\frac{r(\frac{N}{M})}{1-r(\frac{N}{M})}\} + O(\alpha) \right) \right) \nonumber \\
& \gg & M \frac{e^{-u^{2}/2}}{u} \Phi\left(\frac{\alpha/2 - \kappa b (1+O(\alpha))/N^{\alpha} }{\sqrt{(1+O(\alpha))b \alpha/N^{\alpha}}}\right) \min\left\{1,\frac{\alpha N^{\alpha}}{\kappa b}\right\} \sqrt{\frac{\alpha}{b}} \nonumber \\
&& \cdot \prod_{j \leq M^{1/4}} \Phi\left(\left(1+O\left(\frac{\alpha}{j^{\alpha}}\right) \right)\sqrt{\frac{j^{\alpha}}{b\alpha}} \left(1 + \kappa b  \min\left\{\frac{1}{j^{\alpha}},\frac{1}{N^{\alpha}}\right\} + O(\alpha) \right) \right) \nonumber \\
& \gg & M \frac{e^{-u^{2}/2}}{u} \Phi\left(-(1+O(\alpha))\kappa \sqrt{\frac{b}{\alpha Y}} \right) \alpha^{3/2} \prod_{j \leq M^{1/4}} \Phi\left((1+O(\alpha))\sqrt{\frac{j^{\alpha}}{b\alpha}} \left(1 + \kappa b  \min\left\{\frac{1}{j^{\alpha}},\frac{1}{Y}\right\} \right) \right) . \nonumber
\end{eqnarray}
Here the final inequality used the fact that $N^{\alpha}/(\kappa b) = Y/(\kappa b) \gg 1$ and $\sqrt{1/b} \gg 1$, with absolute implied constants, because of our assumptions that $1 \leq b \leq 100$, $1/1000 \leq \kappa \leq 1000$ and $1 \leq Y \leq 1000$.
\begin{flushright}
Q.E.D.
\end{flushright}

\subsection{Proof of Proposition 1}
It remains to prove Proposition 1, and with it Corollary 1, by making a good choice of the remaining parameters $b,\kappa,Y$ in Proposition 2. In order to do this we need to put the lower bound from Proposition 2 into a bit more explicit form.

Whenever $x \geq 2$ we have
$$ \Phi(x) \geq 1 - \frac{1}{x} e^{-x^{2}/2} \geq \exp\{-\frac{2}{x} e^{-x^{2}/2}\}, \;\;\;\;\; \text{and} \;\;\;\;\; \Phi(-x) \gg \frac{1}{x} e^{-x^{2}/2} , $$
and so the lower bound from Proposition 2 implies $\p\left(\max_{1 \leq i \leq M} Z(i/M) > u \right)$ is
\begin{eqnarray}
& \gg & M \frac{e^{-u^{2}/2}}{u} \alpha^{3/2} \Phi\left(-(1+O(\alpha))\kappa \sqrt{\frac{b}{\alpha Y}} \right) \prod_{j \leq Y^{1/\alpha}} \Phi\left((1+O(\alpha))\sqrt{\frac{j^{\alpha}}{b\alpha}} \left(1 + \frac{\kappa b}{Y} \right) \right) \nonumber \\
&& \cdot \prod_{Y^{1/\alpha} < j \leq M^{1/4}} \Phi\left((1+O(\alpha))\sqrt{\frac{j^{\alpha}}{b\alpha}} \left(1 + \frac{\kappa b}{j^{\alpha}} \right) \right) \nonumber \\
& \gg & M \frac{e^{-u^{2}/2}}{u} \alpha^{2} e^{-(1+O(\alpha)) \kappa^{2}b/(2\alpha Y)} \exp\left\{- O\left( \sum_{j \leq Y^{1/\alpha}} \sqrt{\frac{\alpha}{j^{\alpha}}} e^{-(1+O(\alpha))\frac{j^{\alpha}}{2b\alpha}(1+\kappa b/Y)^{2}} \right)\right\} \nonumber \\
&& \cdot \exp\left\{- O\left( \sum_{Y^{1/\alpha} < j \leq M^{1/4}} \sqrt{\frac{\alpha}{j^{\alpha}}} e^{-(1+O(\alpha))\frac{j^{\alpha}}{2b\alpha}(1+\kappa b/j^{\alpha})^{2}} \right)\right\} .\nonumber
\end{eqnarray}
Note that, because $\alpha$ is small, all of the arguments of $\Phi$ had absolute value at least 2, as required.

Next, for any constant $\lambda > 0$ we have that
$$ \sum_{j=1}^{\infty} e^{-\lambda j^{\alpha}/\alpha} \leq \int_{0}^{\infty} e^{-\lambda t^{\alpha}/\alpha} dt = \frac{1}{\lambda} \int_{0}^{\infty} e^{-y} \left(\frac{\alpha y}{\lambda}\right)^{1/\alpha - 1} dy = \frac{1}{\lambda^{1/\alpha}} \alpha^{1/\alpha - 1} \Gamma(1/\alpha) , $$
on substituting $y=\lambda t^{\alpha}/\alpha$. By Stirling's formula the right hand side is $ \ll \frac{1}{\lambda^{1/\alpha}} \alpha^{-1/2} e^{-1/\alpha} $, and so we have
\begin{eqnarray}
\sum_{j \leq Y^{1/\alpha}} \sqrt{\frac{\alpha}{j^{\alpha}}} e^{-(1+O(\alpha))\frac{j^{\alpha}}{2b\alpha}(1+\kappa b/Y)^{2}} & \leq & \sqrt{\alpha} \sum_{j=1}^{\infty} e^{-(1+O(\alpha))\frac{j^{\alpha}}{2b\alpha}(1+\kappa b/Y)^{2}} \nonumber \\
& \ll & \frac{1}{((1+O(\alpha))(1+\kappa b/Y)^{2}/(2b))^{1/\alpha}} e^{-1/\alpha} \nonumber \\
& \ll & \left(\frac{2b}{(1+\kappa b/Y)^{2} e}\right)^{1/\alpha}, \nonumber
\end{eqnarray}
since $(1+O(\alpha))^{1/\alpha} \ll 1$. We remark that it may seem wasteful to remove the factor $1/\sqrt{j^{\alpha}}$, and extend the sum to infinity, but this isn't really the case (given the other parts of our bound): we expect most of the contribution to the sum to come from fairly small $j^{\alpha}$ (because that is where most of the contribution to the gamma function integral arises), and the sum over $Y^{1/\alpha} < j \leq M^{1/4}$ in the other part of our bound is anyway larger than the corresponding part of the sum here.

Next, for any constants $\lambda, \mu > 0$ we have that
$$ \sum_{j=1}^{\infty} \frac{1}{\sqrt{j^{\alpha}}} e^{-(\lambda j^{\alpha} + \mu j^{-\alpha})/\alpha} \ll e^{O(\lambda)} \int_{1}^{\infty} \frac{1}{t^{\alpha/2}} e^{-(\lambda t^{\alpha} + \mu t^{-\alpha})/\alpha} dt , $$
since when $j \leq t \leq j+1$ we see $\lambda t^{\alpha}/\alpha = \lambda j^{\alpha}(1+O(1/j))^{\alpha}/\alpha = \lambda j^{\alpha}/\alpha + O(\lambda j^{\alpha -1})$, and $\mu t^{-\alpha}/\alpha \leq \mu j^{-\alpha}/\alpha$. Substituting $y= t^{\alpha}$, we obtain that
$$ \sum_{j=1}^{\infty} \frac{e^{-(\lambda j^{\alpha} + \mu j^{-\alpha})/\alpha}}{\sqrt{j^{\alpha}}} \ll \frac{e^{O(\lambda)}}{\alpha} \int_{1}^{\infty} \frac{1}{\sqrt{y}} e^{-(\lambda y + \mu y^{-1})/\alpha} y^{1/\alpha - 1}dy = \frac{e^{O(\lambda)}}{\alpha} \int_{1}^{\infty} \left( e^{-(\lambda y + \mu y^{-1})} y\right)^{1/\alpha} \frac{dy}{y^{3/2}} , $$
and calculus shows that the maximum of $e^{-(\lambda y + \mu y^{-1})} y$ over $y > 0$ occurs when $y=(1/2\lambda)(1 + \sqrt{1+4\lambda \mu})$, and is equal to $e^{-\sqrt{1+4\lambda \mu}} (1/2\lambda)(1 + \sqrt{1+4\lambda \mu})$. Therefore\footnote{With more work one could sharpen this bound a bit, by showing that $e^{-(\lambda y + \mu y^{-1})} y$ is only close to its maximum on a small range of $y$. But this wouldn't lead to any real improvement in Proposition 1.}
$$ \sum_{j=1}^{\infty} \frac{1}{\sqrt{j^{\alpha}}} e^{-(\lambda j^{\alpha} + \mu j^{-\alpha})/\alpha} \ll \frac{e^{O(\lambda)}}{\alpha} \left(e^{-\sqrt{1+4\lambda \mu}} (1/2\lambda)(1 + \sqrt{1+4\lambda \mu}) \right)^{1/\alpha} , $$
and so
\begin{eqnarray}
\sum_{Y^{1/\alpha} < j \leq M^{1/4}} \sqrt{\alpha} \frac{e^{-(1+O(\alpha))\frac{j^{\alpha}}{2b\alpha}(1+\kappa b/j^{\alpha})^{2}}}{\sqrt{j^{\alpha}}} & = & \sqrt{\alpha} e^{-(1+O(\alpha))\frac{\kappa}{\alpha}} \sum_{Y^{1/\alpha} < j \leq M^{1/4}} \frac{e^{-(1+O(\alpha))((1/2b)j^{\alpha} + (\kappa^{2}b/2)j^{-\alpha})/\alpha}}{\sqrt{j^{\alpha}}} \nonumber \\
& \ll & \frac{1}{\sqrt{\alpha}} e^{O(\kappa + 1/2b)} \left(e^{-\kappa} e^{-\sqrt{1+\kappa^{2}}}b(1+\sqrt{1+\kappa^{2}}) \right)^{1/\alpha} , \nonumber
\end{eqnarray}
using as before the fact that $(1+O(\alpha))^{1/\alpha} \ll 1$.

Now, for ease of writing, let us temporarily set $F = F(u,\alpha) := \frac{u^{2/\alpha}}{2^{1/\alpha}} \frac{e^{-u^{2}/2}}{u} \alpha^{1/\alpha}$. Then in summary, bearing in mind our restrictions that $1 \leq b \leq 100$, $1/1000 \leq \kappa \leq 1000$ and $1 \leq Y \leq 1000$ (which imply e.g. that $e^{O(\kappa + 1/2b)} \ll 1$), and the fact that $M = \lfloor (bu^{2}\alpha/2)^{1/\alpha} \rfloor$, we have shown that
\begin{eqnarray}
&& \p\left(\max_{1 \leq i \leq M} Z(i/M) > u \right) \nonumber \\
& \gg & M \frac{e^{-u^{2}/2}}{u} \alpha^{2} e^{-\kappa^{2}b/(2\alpha Y)} \exp\left\{- O\left( \left(\frac{2b}{(1+\frac{\kappa b}{Y})^{2} e}\right)^{1/\alpha} + \frac{1}{\sqrt{\alpha}} \left(e^{-\kappa -\sqrt{1+\kappa^{2}}}b(1+\sqrt{1+\kappa^{2}}) \right)^{1/\alpha}  \right)\right\} \nonumber \\
& \gg & F \alpha^{2} \left(b e^{-\kappa^{2}b/2Y}\right)^{1/\alpha} \exp\left\{- O\left( \left(\frac{2b}{(1+\frac{\kappa b}{Y})^{2} e}\right)^{1/\alpha} + \frac{1}{\sqrt{\alpha}} \left(e^{-\kappa-\sqrt{1+\kappa^{2}}}b(1+\sqrt{1+\kappa^{2}}) \right)^{1/\alpha}  \right)\right\} . \nonumber
\end{eqnarray}

\vspace{12pt}
Now we can make our grand selection of parameters. Firstly we must ensure that the second bracket inside the ``big Oh'' term is $< 1$, and this will hold provided $b < f(\kappa)$, where
$$ f(\kappa) := \frac{e^{\kappa + \sqrt{1+\kappa^{2}}}}{1+\sqrt{1+\kappa^{2}}} . $$
Note that if $b$ is strictly smaller than $f(\kappa)$ then the bracketed term will kill off the prefactor $1/\sqrt{\alpha}$. We also remark that $f(\kappa)$ is an increasing function of $\kappa \geq 0$, so we will certainly have $f(\kappa) \geq f(0) = e/2$.

We must also ensure that the first bracket inside the ``big Oh'' term is $\leq 1$, which will hold provided
$$ 1+\frac{\kappa b}{Y} \geq \sqrt{2b/e} . $$
If these two conditions are satisfied (together with the previous restrictions that $1 \leq b \leq 100$, $1/1000 \leq \kappa \leq 1000$ and $1 \leq Y \leq 1000$) then we will have
$$ \p\left(\max_{1 \leq i \leq M} Z(i/M) > u \right) \gg F \alpha^{2} \left(b e^{-\kappa^{2}b/2Y}\right)^{1/\alpha} = \frac{u^{2/\alpha}}{2^{1/\alpha}} \frac{e^{-u^{2}/2}}{u} \alpha^{2} \left(\alpha b e^{-\kappa^{2}b/2Y}\right)^{1/\alpha} . $$

To obtain the best possible lower bound, we should clearly choose $Y$ as large as possible (for given $\kappa$ and $b$), so we choose $Y$ such that
$$ \frac{\kappa b}{Y} = \sqrt{2b/e} - 1. $$
Assuming this choice satisfies $1 \leq Y \leq 1000$, (which it will, for the choices of $\kappa$ and $b$ that we shall make), we will then have
$$ \p\left(\max_{1 \leq i \leq M} Z(i/M) > u \right) \gg \frac{u^{2/\alpha}}{2^{1/\alpha}} \frac{e^{-u^{2}/2}}{u} \alpha^{2} \left(\alpha b e^{-\frac{\kappa}{2} (\sqrt{2b/e}-1)}\right)^{1/\alpha} . $$

For any given $b > e/2$ (which it certainly must be to possibly prove Proposition 1) the above bound is maximised by choosing $\kappa$ as small as possible. And we must always satisfy the constraint $f(\kappa) > b$, so the best lower bound we can possibly obtain is
$$ \frac{u^{2/\alpha}}{2^{1/\alpha}} \frac{e^{-u^{2}/2}}{u} \alpha^{2} \left(\alpha \max_{1/1000 \leq \kappa \leq 1000} f(\kappa) e^{-\frac{\kappa}{2} (\sqrt{2f(\kappa)/e}-1)} \right)^{1/\alpha} . $$
Although the maximum of $f(\kappa) e^{-\frac{\kappa}{2} (\sqrt{2f(\kappa)/e}-1)}$ surely doesn't have a nice closed form, using numerical methods we find it is attained when $\kappa \approx 1.18267$. If we set $\kappa = 1.18267$ then we find $f(1.18267) \approx 6.02449$, so we can choose $b=6.02448$, say. Then if we choose $Y$ such that $\frac{\kappa b}{Y} = \sqrt{2b/e} - 1$, we check that $Y \approx 6.446$, which is also permissible. So, finally, we obtain that
$$ \p\left(\max_{1 \leq i \leq M} Z(i/M) > u \right) \gg \frac{u^{2/\alpha}}{2^{1/\alpha}} \frac{e^{-u^{2}/2}}{u} \alpha^{2} \left(\alpha 6.02448 e^{-0.591335 (\sqrt{12.04896/e}-1)}\right)^{1/\alpha} , $$
and the quantity in brackets is $\approx (3.13362 \alpha)^{1/\alpha} > (1.15279 e \alpha)^{1/\alpha}$. This completes the proof of Proposition 1, and hence of Corollary 1.
\begin{flushright}
Q.E.D.
\end{flushright}

We conclude by making a general remark about the foregoing calculations. Increasing the value of $b$ increases the size of $M$, which increases the number of sample points $i/M$ that are used to lower bound the continuous maximum over $[0,1]$. However, the lower bound supplied by Proposition 2 (ultimately coming from Theorem 2) will deteriorate when $b$ becomes very large, which might be regarded as an undesirable feature of our method. (The cause is that the application of Slepian's lemma in the proof of Theorem 2 becomes increasingly wasteful when $b$ becomes very large). Nevertheless, by optimising $b$ at the same time as optimising the Brownian motion boundary path (i.e. the parameters $C,K,N$) one obtains rather strong lower bounds. The complex form of the bounds as a function of the parameters reflects contributions from different parts of the Brownian motion boundary path, and perhaps also the underlying complexity of understanding $\sup_{0 \leq t \leq 1} Z(t)$.

\appendix

\section{Proofs of the Brownian motion lemmas}
In this appendix we shall prove Brownian Motion Lemmas 1 and 2, as stated in $\S 2.2$. Both proofs exploit a well known explicit formula for the hitting time of a line by Brownian motion, which states that if $\{W_{s}\}_{s \geq 0}$ is a standard Brownian motion started from 0, and if $a > 0$, and if $t > 0$ and $b \in \R$, then
$$ \p(W_{s} \leq a+bs \; \forall 0 \leq s \leq t) = \Phi(\frac{a+bt}{\sqrt{t}}) - e^{-2ab}\Phi(\frac{bt-a}{\sqrt{t}}) , $$
where $\Phi$ denotes the standard normal cumulative distribution function. This formula follows by studying the distribution of the maximum (up to time $t$) of Brownian motion with a drift. See e.g. Chapters 13.4--13.5 of Grimmett and Stirzaker~\cite{gs}.

\subsection{Proof of Brownian Motion Lemma 1}
Let $B:= -b > 0$ (under the hypotheses of the lemma), and let us rewrite the explicit formula as
$$ \p(W_{s} \leq a+bs \; \forall 0 \leq s \leq t) = \frac{1}{\sqrt{2\pi}} \int_{(Bt-a)/\sqrt{t}}^{\infty} e^{-z^{2}/2} dz - e^{2aB} \frac{1}{\sqrt{2\pi}} \int_{(Bt+a)/\sqrt{t}}^{\infty} e^{-z^{2}/2} dz . $$

Suppose first that $a \geq Bt/2$, say. Then we certainly have
\begin{eqnarray}
\int_{(Bt+a)/\sqrt{t}}^{\infty} e^{-z^{2}/2} dz \leq \frac{1}{(Bt+a)/\sqrt{t}} \int_{(Bt+a)/\sqrt{t}}^{\infty} z e^{-z^{2}/2} dz & = & \frac{1}{(Bt+a)/\sqrt{t}} e^{-(Bt+a)^{2}/2t} \nonumber \\
& = & \frac{e^{-2aB}}{(Bt+a)/\sqrt{t}} e^{-(Bt-a)^{2}/2t} \nonumber \\
& \leq & \frac{(2/3) e^{-2aB}}{B\sqrt{t}} e^{-(Bt-a)^{2}/2t} . \nonumber
\end{eqnarray}
On the other hand, it is easy to check that if $a \geq Bt/2$ and $B\sqrt{t}$ is sufficiently large (as hypothesised in Brownian Motion Lemma 1),
$$ \int_{(Bt-a)/\sqrt{t}}^{\infty} e^{-z^{2}/2} dz \geq \frac{1}{B\sqrt{t}} e^{-(Bt-a)^{2}/2t} . $$
Indeed this follows from integration by parts if $(Bt-a)/\sqrt{t} \geq 10$, say, and it is trivial otherwise. So provided $a \geq Bt/2$ and $B\sqrt{t}$ is large we have
$$ \p(W_{s} \leq a+bs \; \forall 0 \leq s \leq t) \geq (1/3)\Phi(\frac{a+bt}{\sqrt{t}}) \gg \min\left\{1,\frac{a}{|bt|}\right\} \Phi(\frac{a+bt}{\sqrt{t}}) , $$
as claimed.

It remains to treat the case where $a < Bt/2$. If we let $\Delta := a/\sqrt{t}$, we note that
$$ \int_{(Bt-a)/\sqrt{t}}^{\infty} e^{-z^{2}/2} dz = \int_{(Bt-a)/\sqrt{t}}^{\infty} e^{-(z+2\Delta)^{2}/2} e^{\Delta(2(z+2\Delta) - 2\Delta)} dz = \int_{(Bt+a)/\sqrt{t}}^{\infty} e^{-w^{2}/2} e^{\Delta(2w - 2\Delta)} dw . $$
But $e^{\Delta(2w-2\Delta)} \geq e^{2aB}$ for all $w \geq (Bt+a)/\sqrt{t}$, and if $w \geq (Bt+a)/\sqrt{t} + 1/(10B\sqrt{t})$ then
$$ e^{\Delta(2w-2\Delta)} \geq e^{2aB} e^{\Delta/(5B\sqrt{t})} = e^{2aB} e^{a/(5Bt)} . $$
We conclude from all these calculations that
\begin{eqnarray}
&& \int_{(Bt-a)/\sqrt{t}}^{\infty} e^{-z^{2}/2} dz - e^{2aB} \int_{(Bt+a)/\sqrt{t}}^{\infty} e^{-z^{2}/2} dz \nonumber \\
& \geq & (e^{a/(5Bt)}-1)e^{2aB} \int_{(Bt+a)/\sqrt{t} + 1/(10B\sqrt{t})}^{\infty} e^{-z^{2}/2} dz \nonumber \\
& \gg & (e^{a/(5Bt)}-1)e^{2aB} \int_{(Bt+a)/\sqrt{t}}^{\infty} e^{-z^{2}/2} dz \gg (e^{a/(5Bt)}-1) \int_{(Bt-a)/\sqrt{t}}^{\infty} e^{-z^{2}/2} dz , \nonumber
\end{eqnarray}
where the penultimate inequality uses the fact that $(Bt+a)/\sqrt{t} < (3/2)B\sqrt{t}$ (say), and the final inequality follows from integration by parts, similarly as in the preceding paragraph. So we have again shown that
$$ \p(W_{s} \leq a+bs \; \forall 0 \leq s \leq t) \gg (e^{a/(5Bt)}-1) \Phi(\frac{a+bt}{\sqrt{t}}) \gg \min\left\{1,\frac{a}{|bt|}\right\} \Phi(\frac{a+bt}{\sqrt{t}}) . $$
\begin{flushright}
Q.E.D.
\end{flushright}

\subsection{Proof of Brownian Motion Lemma 2}
To prove Brownian Motion Lemma 2 we again distinguish two cases, according as $a/\sqrt{t}$ is large enough in terms of $H$, or not. Firstly, if $a/\sqrt{t}$ is large enough then $(a+bt)/\sqrt{t}$ is large and positive (since $|b\sqrt{t}| \leq H$, by hypothesis), $(bt-a)/\sqrt{t}$ is large and negative, and integration shows that
$$ \p(W_{s} \leq a+bs \; \forall 0 \leq s \leq t) = \Phi(\frac{a+bt}{\sqrt{t}}) - e^{-2ab}\Phi(\frac{bt-a}{\sqrt{t}}) = \Phi(\frac{a+bt}{\sqrt{t}}) + O(\frac{1}{|bt-a|/\sqrt{t}}) \gg 1 , $$
as required.

On the other hand, if $a/\sqrt{t}$ is smaller then we again write $B := -b \geq 0$, and by hypothesis we have $(Bt+a)/\sqrt{t} \ll_{H} 1$. Therefore we see, as in the second part of the proof of Brownian Motion Lemma 1 (with $1/(10B\sqrt{t})$ replaced there by 1), that
\begin{eqnarray}
\int_{(Bt-a)/\sqrt{t}}^{\infty} e^{-z^{2}/2} dz - e^{2aB} \int_{(Bt+a)/\sqrt{t}}^{\infty} e^{-z^{2}/2} dz & \geq & (e^{2a/\sqrt{t}}-1)e^{2aB} \int_{(Bt+a)/\sqrt{t} + 1}^{\infty} e^{-z^{2}/2} dz \nonumber \\
& \gg_{H} & (e^{2a/\sqrt{t}}-1) , \nonumber
\end{eqnarray}
and so indeed
$$ \p(W_{s} \leq a+bs \; \forall 0 \leq s \leq t) \gg_{H} (e^{2a/\sqrt{t}} - 1) \geq \frac{a}{\sqrt{t}} . $$
\begin{flushright}
Q.E.D.
\end{flushright}

%%\vspace{12pt}
%%\noindent {\em Acknowledgements.}

\end{document}